\documentclass[12pt]{article}
\usepackage{amsmath,amsthm,amssymb,amscd,}
\newtheorem{thm}{Theorem}[section]
 \newtheorem{cor}[thm]{Corollary}
 \newtheorem{lem}[thm]{Lemma}
 \newtheorem{prop}[thm]{Proposition}
 \theoremstyle{definition}
 \newtheorem{defn}[thm]{Definition}
 \theoremstyle{remark}
 \newtheorem{rem}[thm]{Remark}
 \newtheorem{ex}{Example}
 \numberwithin{equation}{section}
\def\r{\mathbb{R}}
\def\reel{\mathbb{R}}
\def\complex{\mathbb{C}}

\def\rm {\r_+^* \times M}
\def\rn{\r_-\times M}
\def\tm{N}

\title
{Maximal   regularity  and Hardy spaces
}
\author{ Pascal Auscher, Fr\'ed\'eric Bernicot and Jiman Zhao
  \\
 \small Universit\'e de Paris-Sud, Orsay et CNRS 8628, 91405 Orsay Cedex, France
\\
\small {\em E-mail address:} {pascal.auscher@math.u-psud.fr}\\
\small Universit\'e de Paris-Sud, Orsay et CNRS 8628, 91405 Orsay Cedex, France
\\
\small {\em E-mail address:} {Frederic.Bernicot@math.u-psud.fr}\\
\small School of Mathematical Sciences,
Beijing Normal University \\
\small Beijing 100875,
P.R. China\\
\small {\em E-mail address:} {jzhao@bnu.edu.cn}\\ 
}
\date{October 10, 2007}
\begin{document}
\maketitle

\begin{abstract} In this work, we consider  the Cauchy problem for $u' - Au  = f$ with $A$  the Laplacian operator   on some Riemannian manifolds or a sublapacian on some Lie groups or some second order elliptic operators on a domain.  We show the boundedness of the operator of   maximal  regularity  $f\mapsto Au$ and its adjoint on appropriate Hardy spaces which we define and study for this purpose. As a consequence we reobtain the   maximal   $L^q$    regularity on $L^p$ spaces for $1<p,q<\infty$.

\vspace{0.1in}

\vspace{0.1in}
\end{abstract}

\section{Introduction}
 
Let $B$ be a Banach space, $A$ the infinitesimal generator of a bounded analytic semigroup of operators on $B$ and  $I=(0,+\infty)$ ($I$ bounded is also of interest).  We shall consider the problem 
\begin{equation}\label{eq:intro}\left\{\begin{array}{ll}
u'(t)-Au(t)=f(t), &  t \in I,\\
u(0)=0, & 
\end{array}\right. 
\end{equation}
where $f: I\rightarrow B$ is given.
If $T_t$ is the semigroup generated by $A$, then $u$ is 
given by
$$u(t)=\int_0^t T_{t-s}f(s)ds.$$
For fixed $q \in (1,+\infty)$, one says that there is   maximal  $L^q$   regularity
 if for every $f \in
L^q(I, B)$, $u'$ and $Au$
belongs to $L^q(I, B)$ with
\begin{equation}
\label{eq:maxreg}
\|u'\|_{L^q(I,B)} +\|Au\|_{L^q(I,B)} \le C \|f\|_{L^q(I,B)}.
\end{equation}
 When $B$ is an $L^p$ space ($1<p<\infty$) we refer 
to \cite{CV}, \cite{CL}, \cite{CD},  \cite{DS}, \cite{HP}, \cite{L} etc.  When $B$ is an UMD space,  \cite{W} gives a necessary and sufficient condition for    maximal   $L^q$    regularity in terms of a notion called R-boundedness.  An excellent survey is  \cite{kw}.

Let us come back to $B=L^p(\Omega,d\mu)$. A classical method to prove maximal regularity (already in \cite{DS}) relies on singular integral theory.  It is well-known that  maximal   $L^q$    regularity for one $q\in (1,+\infty)$ implies the same for all $q\in (1,+\infty)$.  Hence, it is enough to prove \eqref{eq:maxreg} for $q=p$, that is the boundedness of the map of  maximal regularity $f\mapsto Au$ on $L^p(I\times \Omega, dtd\mu)$. Seen as a singular integral on $I\times \Omega$ equipped with a parabolic distance and product measure $dtd\mu$, one proves this map and its adjoint have weak type (1,1)  and applies interpolation. This method has been refined so as to relax hypotheses as much as possible on $\Omega$ and $A$. Completing a theorem in \cite{HP}, this method has been successful  for all $1<p<\infty$ in
 \cite{CD} with $\Omega$ a (subset of a) space of homogeneous type and pointwise Gaussian upper bounds for the kernel of the semigroup and even generalized to a restricted range of $p$ in \cite{bk} under weaker generalized  Gaussian estimates. 
 
 There exists a criterion for $L^1$ maximal regularity on $L^1$  (see \cite{gd}) but when $A$ is a second order differential operator it does not apply. Of course, seen from the point of view of singular integrals, the map $f\mapsto Au$ is naturally not bounded on $L^1(I\times \Omega)$ and one should replace $L^1(I\times \Omega)$ by a Hardy space and prove  boundedness  into $L^1(I \times \Omega)$. (The notion of $H^p$-maximal regularity  introduced in \cite{blm} is a different problem.) When $A$ is the Laplacian on $\r^n$, kernel representation suggests to use a parabolic Hardy space of Coifman-Weiss type on $I\times  \r^n$. Indeed, this is done implicitly in \cite[Appendix]{SV} given $H^1-BMO$ duality,  where the authors prove boundedness from $L^\infty(\r \times \r^n)$ to  parabolic $BMO$  of a dual problem.  They also do this for some time dependent generalization of \eqref{eq:intro} with the Laplacian on $\r^n$. 
 
 In the abstract setting of \cite{CD} and \cite{bk},  it is not clear which Hardy space is  appropriate \footnote{This will be done in a subsequent paper by the last two authors.}. Here, we restrict ourselves to  
 the Laplace-Beltrami operator $\Delta$ on some manifolds $M$ or a sublaplacian on some Lie groups, or second order operators on $\r^n$ with bounded measurable coefficients with additional assumptions.  We prove that the operator of maximal regularity is not just bounded from a Hardy space into $L^1$ but  into itself. We prove a similar phenomenon for the dual operator 
 but with two different Hardy spaces as source and target.
 The Hardy spaces  on which we prove boundedness are of course of parabolic type. Interpolation applies to reobtain  maximal  $L^q$   regularity on $L^p$ for $1<p,q<\infty$.

The plan of the article is as follows. We give next the notation and then state the main result in the setting of a Riemannian manifold. We define and study the parabolic Hardy spaces  in Section 3. Then the proof of the main result on manifolds  is presented  in Section 4 and 5. The consequence for the  maximal   $L^q$   regularity is in Section 6.  We give in Section 7 the needed ingredients to prove the analogous results in a connected Lie group with polynomial growth. We also present the analogous results for a class of second order operators on domains of $\r^n$ in Section 8.


\section{Main result on manifolds}

Let $M$ be a complete non-compact connected Riemannian manifold,
$\mu$ the Riemannian measure. Denote by
$\|\cdot\|_{p}$ the norm in $L^p(M,\mu),~1\leq p \leq
+\infty$. Denote by $\nabla$ the Riemannian
gradient, $|\cdot|$ the length in the tangent space, $\Delta$ the (negative)
 Laplace-Beltrami operator corresponding to $\nabla$ and $d$  the
 geodesic distance. For all $x \in M,$ all $r>0$, $B(x,r)$ stands for
 the open geodesic ball with center $x$ and radius $r$, and its
 measure is  denoted $V(x,r)$.

Say that $M$ satisfies the doubling property if
and only if there exists $C > 0$ such that for all $r>0$ and $x\in M$,
$$V(x,2r) \leq C V(x,r). \eqno(D)$$
A straightforward consequence of $(D)$ is that there exist $C,\delta >0$
such that for all $x \in M$, all $r>0$ and all $\theta>1$,
$$V(x,\theta r) \leq C \theta^\delta  V(x,r).\eqno(D')$$
The hypothesis $(D)$ exactly means that $M$, equipped with its
geodesic distance and its Riemannian measure, is a space of
homogeneous type in the sense of Coifman and Weiss (\cite{cw}).

Say that $M$ satisfies the (scaled) Poincar\'e
inequalities if there exists $C>0$ such that for every ball
$B=B(x,r)$ and every $f$ with $f$, $\nabla f$ locally square
integrable,
$$\int_B|f-f_B|^2 d\mu \leq C r^2\int_B|\nabla
f|^2 d\mu,\eqno(P)$$ where $f_B$ denotes the mean of $f$ on $B$.

The assumptions  $(D)$ and $(P)$ guarantee the needed properties on the
heat kernel (size and    regularity) in the sequel.

Consider the Cauchy problem
\begin{equation}\label{eq0}
\left\{\begin{array}{ll}
\frac{\partial u}{\partial t}(t,x)-\Delta u(t,x)=f(t,x), &  t>0,\\
u(0,x)=0. & 
\end{array}\right.
\end{equation}


Now define 
\begin{equation}\label{eq:T}Tf(t,x)=\int_0^t \big[\Delta
e^{(t-s)\Delta} f(s,\cdot)\big](x)\, ds. 
\end{equation} 
By \cite{DS}, we know that $T$ is bounded from $L^2(X)$ to $L^2(X),$ where $X=\rm$ and  product  measure $dtd\mu$.
Its adjoint is given by
\begin{equation}\label{eq:T*}T^*f(t,x)=\int_t^\infty \big[\Delta
e^{(s-t)\Delta} f(s,\cdot)\big](x)\,ds.\end{equation}

We introduce in the next section two Hardy spaces of parabolic type on $X$ with 
$H^1_z(X)\subset H^1_r(X) \subset L^1(X)$. Our main result is the following.

\begin{thm}\label{thm:h1reg}  Assume that $M$ satisfies $(D)$ and $(P)$. Then 
 $T$ is $H_z^1(X)$-bounded and $T^*$ is bounded from $H_r^1(X)$ to $H_z^1(X)$.
\end{thm}

\begin{cor} With the hypotheses above, one has
$$\|u'\|_{H^1_z(X)} + \|\Delta u \|_{H^1_z(X)} \le  C \|f\|_{H^1_z(X)}.$$
\end{cor}

 \noindent {\bf Proof:}  If $f \in H^1_z(X) \cap L^2(X)$, we have $\Delta u =Tf$ and $u'=f+\Delta u$ by uniqueness of the Cauchy problem. Thus this estimate holds for those $f$'s and we conclude by density in $H^1_z(X)$ (See Section \ref{sec:hardy}).
 \qed
 
 
 We shall also show that a better result cannot hold even in the simplest situation $X=\r_+^*\times \r^n$: $T$ is not bounded from $ H^1_r(X)$ to $ L^1(X)$ and  $T^*$ is not bounded on $L^1(X)$.

\section{Parabolic  $H_r^1(X)$ and $H_z^1(X)$}\label{sec:hardy}

In this section we study the parabolic Hardy spaces $H_r^1(X)$ and
$H_z^1(X)$, following the ideas of \cite{cks} and \cite{ar} in an elliptic context.

We assume that $(E,d,\mu)$ is a space of homogeneous type.
Let $\tm$  be the product space $\r \times E$ endowed with the measure $\nu$,
product of the Lebesgue measure with $\mu$, and with the
quasi-distance $\widetilde d$ defined  by
$$\widetilde{d}(t,x;s,y)^2=\sup (d(x,y)^2, |t-s|), \qquad (t,x),(s,y) \in \tm.$$ 

Suppose $Q =\{ (t,x) \in \tm: \widetilde{d}(t,x;s,y)< r\}$ is the ball in $\tm$ centered at $(s,y)$ with radius $r$. Then the relation between the volumes of $Q$ and $B(y, r)$ is  $\nu(Q) \sim r^2 V(y,r),$ where $V(y,r)=\mu(B(y,r))$.  One then easily  checks  that $(\tm, \widetilde{d}, \nu)$ is a space of homogeneous type and 
\begin{equation}\label{eq:doublingnu}
\nu(\theta Q)\le C\theta^{\delta+2}\nu(Q)
\end{equation}
 for all balls $Q$ and all $\theta>1$ where $\delta$ is the doubling exponent for $E$ as in (D'). 
Thus, one can consider the Hardy space $H^1(\tm)=H^1_{CW}(\tm)$ of \cite{cw}. Let us review some definitions and facts.

\begin{defn}
A classical $(1,\infty)$-atom on $\tm $ is a measurable function
$a$ supported on a ball $Q \subseteq \tm $ that
satisfies

(1) $\int_N ad\nu = 0,$

(2) $\|a\|_{L^\infty(\tm)} \leq
\nu(Q)^{-1}.$
\end{defn}

\begin{defn}
An $(1,2)$-atom on $\tm$ is a measurable function
$a$ supported on a ball $Q \subseteq\tm$ that
satisfies

(1) $\int_N ad\nu=0,$

(2) $\|a\|_{L^2(\tm)} \leq
\nu(Q)^{-\frac{1}{2}}.$

\end{defn}

An $L^1(\tm)$ function $f$ belongs to $H^1(\tm)$ if it can be written as $\sum \lambda_Q a_Q$
where $\sum |\lambda_Q|<\infty$ and $a_Q$ are $(1, 2)$-atoms. The norm is given by the infimum of such $\sum |\lambda_Q|$ taken over all possible decompositions called atomic decompositions. One obtains the same space with an equivalent norm by replacing $(1,2)$-atoms by $(1,\infty)$-atoms.

We now introduce two Hardy spaces on $X=\r_+^*\times E$. 

\begin{defn} Set
$H_r^1(X)=\{f \in L^1(X):  \exists \,F \in H^1(\tm) \, , \,  F|_{X} =f\} .$
The norm of an element  in this space is the quotient norm : the
infimum of the $H^1(\tm)$ norms of all possible extensions  to
$\tm$.
\end{defn}

\begin{defn} Set $H_z^1(X) = \{f|_X : f  \in H^1(\tm), f=0~~\mbox{on}~~\r_- \times E \}.$
The norm of an element  in this space is the norm of its zero extension in $H^1(\tm)$.
\end{defn}

It is clear that $H_z^1(X)$ $\subset$ $H_r^1(X)$, but the
converse is not true as shown by the following example.

\begin{ex} Let
$$
f(t,x)=\left\{\begin{array}{ll}
1,&  (t,x) \in (0,1) \times B(0,1),\\
0,& \mbox{elsewhere},
\end{array}\right.
$$
then $f$ is the restriction to $X$ of
$$
F(t,x)=\left\{\begin{array}{ll}

1,&  (t,x) \in (0,1)\times B(0,1),\\

-1,&  (t,x) \in (-1,0) \times B(0,1),\\

0,& \mbox{elsewhere}.

\end{array}\right.
$$
It is clear that $F$ is a $(1,\infty)-$atom on $\tm$, so $f \in
H_r^1(X)$, but $f \notin H_z^1(X)$ because $f$ does not have vanishing 
moment.
\end{ex}


%


\begin{defn}
Let $a$ be a measurable function supported in a ball $Q \subset
X$, and satisfying $\|a\|_{L^2(\tm)} \leq
\nu(Q)^{-\frac{1}{2}}$.

(1) $a$ is a type $(a)$ atom if, moreover,  $4Q \subseteq X$ and
$\int_X a d\nu =0$.

(2) $a$ is a type $(b)$ atom if, moreover,  $2Q \subseteq X$ and  $4Q \nsubseteq X$ ($a$ does not necessarily have vanishing moment ).
\end{defn}
We have an atomic characterization of $H_r^1(X).$ The proof, analogous to the one in \cite{cks}, is  included  for completeness.
\begin{thm}
Let $f$ be an integrable function on $X$. The followings are equivalent:

(1) $f \in H_r^1(X)$;

(2) $f$ has an atomic decomposition 

\begin{equation}\label{eqh1r}f=\sum_{\mbox{type (a) atoms}} \lambda_Q a_Q + \sum_{\mbox{type (b) atoms}}
\mu_Q b_Q\end{equation}
with
\begin{equation}\label{eqnormh1r}\sum_{\mbox{type (a) atoms}} |\lambda_Q |+ \sum_{\mbox{type (b) atoms}}
|\mu_Q | < \infty.\end{equation}
Further, the norm on $H^1_r(X)$ is equivalent to the infimum of (\ref{eqnormh1r}) taken over all decompositions (\ref{eqh1r}).
\end{thm}

\noindent {\bf Proof:}  (2)$\Rightarrow$(1). Let $f$ satisfy  \eqref{eqh1r} and \eqref{eqnormh1r}. As $a_Q$ is a type ($a$) atom, then $A_Q$ defined by
$$
A_Q(t,x) = \left\{\begin{array}{ll}
a_Q(t,x),& (t,x) \in X,\\
0,& (t,x) \notin X,\\
\end{array}\right.
$$
 is a (1,2)-atom of $H^1(\tm)$. By definition we conclude  that $a_Q \in H_r^1(X)$.

As $b_Q$ is a type ($b$) atom, then define
$$
B_{\widetilde{Q}}(t,x)=\left\{\begin{array}{ll}
b_Q(t,x),& (t,x) \in Q ,\\
-b_Q(-t,x),& (t,x) \in Q_-,
\end{array}\right.
$$
where $Q_-$ is the reflection of $Q$ across $\{0\} \times M$ given by $(t,x) \mapsto (-t,x)$, and
$\widetilde{Q}$ is a ball such that $Q \bigcup Q_-
\subseteq \widetilde{Q}$ and $\nu(\widetilde Q ) \le C \nu(Q)$. For instance, $Q$ is centered at $(s,y)$ with radius $r$, the ball  $\widetilde Q$  centered at $(0,y)$ with radius $5r$ works. 

Then we have
\[
\begin{array}{l}
\int_{\tm} B_{\widetilde{Q}}d\nu
=\iint_X B_{\widetilde{Q}}(t,x) dt d\mu(x) + \iint_{\rn} B_{\widetilde{Q}}(t,x) dt d\mu(x)\\[0.1in]
=\iint_Q b_Q(t,x) dt d\mu(x) - \iint_Q b_Q(t,x) dt d\mu(x)
=0,\\
\end{array}
\]
and
$$
\|B_{\widetilde{Q}}\|_{L^2(\tm)} \leq \|b_Q \|_{L^2(\tm)} + \|-b_Q \|_{L^2(\tm)}
\leq \nu(Q)^{-\frac{1}{2}} + \nu(Q_-)^{-\frac{1}{2}} \leq
c\nu(\widetilde{Q})^{-\frac{1}{2}}.
$$
This means that $c^{-1}B_{\widetilde{Q}}$ is an $(1,2)$-atom for $H^1(\tm)$. Now let 
$$F=\sum_{\mbox{type (a) atoms}} \lambda_Q A_Q + \sum_{\mbox{type (b) atoms}}
\mu_Q B_{\widetilde{Q}},$$
then $F \in H^1(\tm)$ and $F|_X = f$, therefore $f \in H^1_r(X)$.\\

(1) $\Rightarrow$ (2). If $f \in H^1_r(X)$, then there exists $F \in H^1(\tm)$ such that $F|_X=f$. 
So $F$ has an atomic decomposition $\sum \lambda_Q A_Q$ \cite{cw}. It suffices therefore to concentrate on the restriction of one single atom $A_Q$. Restricting $A_Q$ to $X,$ then we consider only those balls $Q$ which intersect $X$ 
 and we have the following situations:

If $4Q \subseteq X$, then $a_Q=A_Q$ is a type ($a$)-atom.

If $4Q \nsubseteq X $, but $2Q  \subseteq X$, then $b_Q=A_Q$
is a  type ($b$)-atom.

If $2Q  \nsubseteq X$, then we can decompose $A_Q|_X$ into
type ($b$)-atoms as follows. 
Using Whitney decomposition on $X$  with respect
to $\r_- \times E$ (see \cite{fs}),  one can find a family of balls $Q_j\subset X$   with $Q\cap X \subseteq \cup Q_j$
 and  bounded overlap. Hence, $\sum \nu( Q_j) \le c\nu(Q\cap X)$ for some constant $c>0$ depending on the metric and the doubling property of $\nu$. Write
$$A=\sum_j \chi_{Q_j} A = \sum_j \frac{\|\chi_{Q_j}A\|_2} {\nu(Q_j)^{-\frac{1}{2}}}
\cdot \frac{\nu(Q_j)^{-\frac{1}{2}}}{\|\chi_{Q_j}A\|_2} \chi_{Q_j}A = \sum
\lambda_{Q_j} a_{Q_j},$$
where the indices $j$ are those for which $\chi_{Q_j}A$ are not identically $0$ and
$$\lambda_{Q_j} = \frac{\|\chi_{Q_j}A\|_2}{\nu(Q_j)^{-\frac{1}{2}}}, ~~a_{Q_j} = \frac{\nu(Q_j)^{-\frac{1}{2}}}{\|\chi_{Q_j}A\|_2}\chi_{Q_j}A.$$
Then we have supp$\,a_{Q_j} \subseteq Q_j$ and 
$$\|a_{Q_j}\|_{L^2(\tm)} \leq \frac{\nu(Q_j)^{-\frac{1}{2}}}{\|\chi_{Q_j}A\|_2}
\cdot {\|\chi_{Q_j}A \|_2} \leq \nu(Q_j)^{-\frac{1}{2}}.$$
It follows that $a_{Q_j}$ is a  type ($b$)-atom and
$$
\sum |\lambda_{Q_j}| = \sum_j \frac {\|\chi_{Q_j}A\|_2}{\nu(Q_j)^{-\frac{1}{2}}}
\leq (\sum_j \|\chi_{Q_j} A\|_2^2)^{\frac{1}{2}} \cdot (\sum_j \nu(Q_j))^{\frac{1}{2}}
\lesssim \|A\|_{L^2(\tm)}\nu(Q\cap X)^{\frac{1}{2}}.
$$
As $\|A\|_{L^2(\tm)} \le \nu(Q)^{-\frac 1 2}$, the series is bounded independently of $A$ and this atomic decomposition converges boundedly in $H^1_r(X)$. \hfill$\Box$

For  the space $H^1_z(X)$, we have the following characterization, whose proof is again inspired by \cite{cks}.

\begin{thm}\label{thm:h1z}
Let $f$ be an integrable function on $X$. The followings are equivalent:

(1) $f \in H^1_z(X)$;

(2) $f$ has an atomic decomposition $f= \sum_{Q \subset X}
\lambda_Q a_Q$ with $\sum |\lambda_Q|< +\infty$,
where $a_Q$ are  (1,2)-atoms with support entirely contained in $X$.
\end{thm}

 \noindent {\bf Proof:} (2) $\Rightarrow$ (1) is obvious and we turn to (1) $\Rightarrow$ (2). Let $f \in H^1_z(X)$. Then
the function  $F$ defined by
$$
F(t,x) = \left\{\begin{array}{ll}
f(t,x),& (t,x) \in X, \\
0,& (t,x) \notin X.
\end{array}\right.
$$
belongs to  $ H^1(\tm)$.
Set $f_e(t,x)=F(t,x)+F(-t,x)$ with $(t,x)\in \tm$
 then
$f_e $ is even in the $t$ variable, $f_e \in H^1(\tm),$ 
$f_e|_X = f$. Pick an
atomic decomposition $f_e = \sum_Q \lambda_Q A_Q$, where
$A_Q$ are $(1,2)$-atoms for $H^1(\tm)$.
One can rewrite $f_e$ as follows:
$$f_e(t,x) = \frac{f_e(t,x) + f_e(-t,x)}{2} = \sum_Q \lambda_Q
\frac {A_Q(t,x) + A_Q(-t,x)}{2} = \sum_Q
\lambda_Q \widetilde{A}_Q(t,x),$$
where $\widetilde{A}_Q(t,x) =\frac{A_Q(t,x) + A_Q(-t,x)}{2}.$

$1^\circ$ If supp$\, A_Q \subseteq X$, then $\widetilde{A}_Q|_X
 = \frac {A_Q}{2}$, and we define
$$a_Q(t,x)  = \frac{A_Q(t,x)}{2}, \quad (t,x) \in X.$$

$2^\circ$ If supp$\,A_Q \subseteq \rn$, $\widetilde{A}_Q|_X(t,x)
 = \frac {A_Q(-t,x)}{2}$, and we define
$$a_Q(t,x)  = \frac{A_Q(-t,x)}{2}, \quad (t,x) \in X.$$

$3^\circ$ If neither case occurs, then let
$Q_0 = $ (supp$\,A_Q) \setminus  X$, and the "reflection" of $Q_0$
is defined by $Q_r =\{ (t,x) \in X:(-t,x) \in Q_0 \}.$ Define $a_Q = \widetilde{A}_Q|_X$.
Clearly, if $Q$ is centered at $(s,y)$ has radius $r$, we have $|s| \le r^2$ and supp$\, a_Q \subset Q|_X \cup Q_r $. Hence supp$\, a_Q $ is contained in the ball $\widetilde Q$ centered at $(r^2,y)$ with radius $r$ which is contained in $X$. Next, 
$$\|a_Q \|_{L^2(\tm)} \leq \nu(Q)^{-\frac{1}{2}}
= \nu(\widetilde Q)^{-\frac{1}{2}}.$$

Finally 
\begin{align*}
\int_X a_Qd\nu
&=\frac{1}{2} \iint_{Q \cap X}  A_Q(t,x) dtd \mu(x)+ \frac{1}{2}\iint_{Q_r }A_Q(-t,x)dtd \mu (x)\\
&=\frac{1}{2} \iint_{Q \cap X} A_Q(t,x)dt d \mu(x) + \frac{1}{2}\iint_{Q \setminus X }A_Q(t,x)dt d \mu (x)\\
&=\iint_Q A_Q(t,x)dtd \mu (x) =0.
\end{align*}
Thus $a_Q$ are $(1,2)$-atoms with supp$\,a_Q \subset X.$
So
$$f = f_e|_X = \sum_Q \lambda_Q \widetilde{A}_Q|_X
=\sum_Q \lambda_Q a_Q.$$
This completes the proof of the theorem. \hfill$\Box$

Let us make some further remarks. If one uses the reflection $(t,x) \mapsto (-t,x)$ in $\tm$, then we say that a function on $\tm$ is odd (resp. even) if it is changed to its opposite (resp. itself).

\begin{prop}\label{prop:reflexion} $ H^1_r(X)$ coincides with the space of restrictions to $X$ of odd functions in $ H^1(\tm)$. Similarly, $H^1_z(X)$ coincides with the space of restrictions to $X$ of even functions in $H^1(\tm)$. 
\end{prop}

\noindent {\bf Proof:}  The argument is essentially contained in the proofs of the last two results and is skipped (See also \cite{cks}). \hfill$\Box$

\begin{rem} Density of $H^1_z(X) \cap L^2(X)$ (resp. $H^1_r(X) \cap L^2(X)$) in $H^1_z(X)$ (resp. $H^1_r(X)$) follows from 
the atomic decomposition.
\end{rem}

Since $(X,\widetilde d,\nu)$ is also a space of homogeneous type, let  
$H^1(X)=H^1_{CW}(X)$ be the Hardy space of Coifman and Weiss as defined in \cite{cw}. Each $f \in H^1(X)$ can be written as follows:
$$f = \sum_Q \lambda_Q a_Q,$$ 
where supp$\,a_Q \subset Q
\cap X$, $Q$ is a ball centered in $X$, $\int _{Q\cap X} a d\nu=0,$ and 
$\|a_Q\|_{L^2(\tm)} \leq {\nu (Q \cap
X)^{-\frac{1}{2}}}$ and where $\sum |\lambda_Q|<\infty$. Observe that since $Q$ is centered in $X$, one has $\nu(Q)\lesssim \nu (Q \cap
X)$ with implicit constant independent of $Q$. 

\begin{prop}\label{prop:h1z=h1cw}
$H^1(X)=H^1_z(X).$
\end{prop}

\noindent {\bf Proof:} That $H^1_z(X) \subseteq H^1(X)$ is obvious.
To prove the converse inclusion let 
$a$ be an  atom in  $ H^1(X)$: there exists a ball $Q$ in $\tm$ centered at a point in $X$ such that   supp$\, a \subseteq Q \cap X$.
If $ Q  \subset X$, then $a \in H^1_z(X)$. Otherwise, let $(s,y)$ be the center of $Q$ and $r$ its radius and observe that  $0<s$. Let $\widetilde Q$ be the ball centered at $(r^2,y)$  with radius $r$.  Then 
$$\nu(\widetilde Q) \sim r^2V(y,r)= \frac{r^2}{s+r^2} (s+r^2)V(y,r) \le 2 \nu(Q\cap X).
$$
Clearly supp$\, a \subseteq \widetilde Q$ and $ a$ is up to a fixed multiplicative constant, an $(1,2)$-atom supported inside $X$. This proves that $a\in H^1_z(X)$ and completes the proof of the proposition. \hfill$\Box$

\begin{rem} One obtains all the above atomic decompositions replacing atoms with $L^2$ estimates by atoms with $L^\infty$ estimates. This will be used in the sequel.
\end{rem}

We now recall the notion of molecules.

\begin{defn}\label{defn:molecule} Let $Q$ be a ball contained in $X$.\\
A function $m$ is called a molecule associated to $Q$ if there
exists $\alpha >0$ such that

(1) $\int_X m  d\nu =0$,

(2) for all $j\ge 1$, $\|m\|_{L^2(B_j(Q))} \leq {2^{-j\alpha}} \nu (2^{j+1}Q \cap X)^{-\frac{1}{2}}$,

where
$$
B_j(Q)=\left\{\begin{array}{ll}
4Q \cap X,&  j=1, \\
(2^{j+1}Q \setminus 2^jQ) \cap X,&  j \geq 2.
\end{array}\right.
$$
A function $m$ is called a  type $(a)$ molecule associated to $Q$ if  $4Q \subseteq X$ and there
exists $\alpha >0$ such that (1) and  (2) hold.\\
A function $m$ is called a  type $(b)$ molecule associated to $Q$ if  $2Q \subseteq X$, $4Q \nsubseteq X$ and there
exists $\alpha >0$ such that   (2) holds.

\end{defn}

\begin{prop}\label{lem:molecule} Molecules form a bounded set in $H^1_z(X)$. Molecules of type (a) and (b) form a bounded set in $H^1_r(X)$.
\end{prop}

Proofs are somewhat  analogous to the ones in \cite{cw} but we include them  for convenience. \

\noindent {\bf Proof:} Let  us begin with the case of $H^1_z(X)$. Let $m$ be a molecule associated to a cube $Q$. By  Proposition \ref{prop:h1z=h1cw}, it suffices to show that $m$ belongs to $H^1(X)$ with norm independent of $Q$.  Set $\chi_{A}$ the indicator of a set $A$ and  $m_{B_j(Q)}$  the mean value of $m$ over $B_j(Q)$: $ \frac{1}{\nu (B_j(Q))}\int_{B_j(Q)}m\, d\nu$. Write
$$
\ m = \sum_{j=1}^\infty m\chi_{B_j(Q)}
= \sum_{j=1}^\infty (m-m_{B_j(Q)})\chi_{B_j(Q)} + \sum_{j=1}^\infty m_{B_j(Q)}\ \chi_{B_j(Q)}.
$$
Using (2) in Definition \ref{defn:molecule}, we obtain that
$$
\int_{2^{j+1}Q\cap{X}}
|(m-m_{B_j(Q)})\chi_{B_j(Q)})|^2d\nu 
\leq 4 \int_{B_j(Q) } |m|^2d\nu
\lesssim \frac{1}{4^{j\alpha}} \cdot \frac{1}{\nu (2^{j+1}Q\cap{X})},
$$
where we have used the doubling property $\nu(2^{j+1}Q) \sim \nu(2^jQ)$ and also the fact that $\nu(Q\cap X) \sim \nu(Q)$ for all balls centered at a point in $X$. 

Set
$$\lambda_j = \frac{1}{2^{j\alpha}}, ~~a_j = {\lambda_j}^{-1}(m - m_{B_j(Q)})\chi_{B_j(Q)}.$$
We have 
$\int_{2^{j+1}Q\cap{X}} a_j d\nu = 0, \mbox{supp}\,a_j \subset 2^{j+1}Q \cap{X},$
and $\|a_j\|_{L^2(2^{j+1}Q \cap{X})} \leq C\nu ( 2^{j+1}Q \cap{X})^{-\frac{1}{2}},$ with $C$ independent of $j$ and $Q$. 
So $a_j \in {H^1}(X)$ and $\|a_j\|_{H^1(X)}$ is bounded. 
Since
$\sum_{j=1}^\infty \lambda_j \le C(\alpha)<\infty,$  
 $$\sum_{j=1}^\infty \left(m- m_{B_j(Q)}\right)\chi_{B_j(Q)}=
\sum_{j=1}^\infty \lambda_j a_j \in  H^1(X).
$$
For the second term, we have 
\[
\begin{array}{l}
\sum_{j=1}^\infty m_{B_j(Q)}\ \chi_{B_j(Q)}=
\displaystyle \sum_{j=1}^\infty \int_{B_j(Q)}m\, d\nu \ \frac{\chi_{B_j(Q)}}{\nu (B_j(Q))}\\[0.1in]
=\displaystyle  \sum_{j=2}^\infty \left(\int_{2^{j+1}Q\cap{X}}m \, d\nu-
\int_{2^jQ\cap{X}}m\, d\nu \right) \frac{\chi_{B_j(Q)}}{\nu( B_j(Q))} +  \left(\int_{4Q\cap{X}}m  \, d\nu \right) \frac{\chi_{B_1(Q)}}{\nu( B_1(Q))} \\[0.1in]
= \displaystyle  \sum_{j=2}^\infty \left( \int_{2^jQ\cap{X}}m\, d\nu \right)
\bigg[\frac{\chi_{B_{j-1}(Q)}}{\nu (B_{j-1}(Q))} - \frac{\chi_{B_j(Q)}}{\nu (B_j(Q))}\bigg].
\end{array}
\]
Set 
$$
\mu_j = \int_{2^jQ\cap{X}} m d\nu,~~
b_j = \frac{\chi_{B_{j-1}(Q)}}{\nu (B_{j-1}(Q))} - \frac{\chi_{B_j(Q)}}{\nu (B_j(Q))}.
$$
Then by $(1)$ in Definition \ref{defn:molecule},
$$\mu_j = \int_{2^jQ\cap{X}} m \, d\nu
=-\left(\int_X m \, d\nu - \int_{2^jQ \cap{X}} m \, d\nu \right)
=-\sum_{i \geq j} \int_{B_i(Q)} m \, d\nu,
$$
so
$$ |\mu_j| \leq \sum_{i\geq j}\int_{B_i(Q)}|m|\, d\nu\lesssim \sum_{i \geq j}{2^{-i \alpha}} \sim  {2^{-j\alpha}},
$$ 
thus
$$\sum_{j=1}^{\infty} |\mu_j| \lesssim
\sum_{j=1}^{\infty} {2^{-j\alpha}} < + \infty.
$$
On the other hand, supp$\, b_j \subseteq 2^{j+1}Q\cap X$, 
$$
\int_{2^{j+1}Q\cap{X}} b_j d\nu =0,$$
and
$$
\int _{2^{j+1}Q\cap{X}} |b_j|^2 d\nu \leq \frac{C^2}{\nu(2^{j+1}Q\cap X)},$$ for some $C>0$ independent of $j$ and $Q$.  To prove the last inequality, we remark that 
$\nu(B_j(Q))$  is comparable to $\nu(2^{j+1}Q\cap X) $ for all $j\ge 1$.  This is obvious for $j=1$ and we turn to $j\ge 2$.  As  
$B_j(Q)$ is contained in $2^{j+1}Q\cap X$, we obtain one inequality. Next, write $Q=I(t_0, r^2) \times B(x_0, r)=I \times B$ and remark that $B_j(Q)$ contains the set  $[t_0+(2^{j}r)^2, t_0 +(2^{j+1}r)^2) \times (2^{j+1}B)$.   Hence
$$
\nu(B_j(Q)) \ge (2^{j}r)^2 \mu(2^{j+1}B) \sim \nu(2^{j+1}Q) \sim \nu(2^{j+1}Q\cap X)
$$
where the last equivalence comes from the fact that $2^{j+1}Q$ is centered in $X$. 

Thus $C^{-1} b_j$ are atoms for $H^1(X)$, and therefore $\sum_{j=1}^\infty \mu_j b_j \in H^1(X).$ 
 
 \
 
  Clearly the same argument applies for type $(a)$ molecules  so they belong to $H^1_z(X) \subset H^1_r(X)$. It remains to consider a type $(b)$ molecule $m$ associated to a cube $Q$. According to Proposition \ref{prop:reflexion}, the odd extension $m_{odd}$ is, up to a constant that depends only on $\tm$, a molecule for $H^1(\tm)$ associated to a cube $\widetilde Q$ of size comparable to $Q$ (since $4Q \nsubseteq X$), containing $Q$ and its reflection across $\{0\} \times E$. The same argument as above shows that $m_{odd}$ belongs to $H^1(\tm)$. Thus $m\in H^1_r(X)$ and the control of its norm follows from an examination of the argument.  \hfill{$\Box$} 
  
  \
  
We finish with a remark on how to prove boundedness on Hardy spaces for an operator, bypassing any other knowledge such as weak type (1,1).
 
\begin{prop}\label{prop:bdd}
Let $H_1$, $H_2$ be any  Hardy spaces obtained via an atomic decomposition. Any linear operator mapping the atoms of $H_1$ into a bounded set in $H_2$ has a bounded extension
from $H_1$ into $H_2$.
\end{prop}

\noindent {\bf Proof:} Let $S$ be the subset of $H_1$ consisting of finite linear combinations 
$f=\sum_{i=1}^n \lambda_i a_i$, $n\ge 1$, $\lambda_i \in \complex$, $a_i$
 atoms of $H_1$ with $\|f\|_{H_1} \ge \frac 1 {10} \sum_{i=1}^n |\lambda_i|$.
 We claim that $S$ is dense in $H_1$. Indeed, if $f\in H_1$, there is a (possibly) infinite representation $f=\sum_{i=1}^\infty  \lambda_i a_i$ with 
 $ \sum_{i=1}^\infty |\lambda_i| \le 2\|f\|_{H_1}$. As the partial sums $f_n$ converge to $f$  in $H_1$, we have $\|f\|_{H_1} \le 5 \|f_n\|_{H_1}$ for $n$ large enough and 
 $$
 \sum_{i=1}^n |\lambda_i| \le 2 \|f\|_{H^1} \le 10   \|f_n\|_{H_1} .
 $$
 This proves that $f_n \in S$ for all $n$ large enough.
 
 By assumption, we have a linear operator $T$ satisfying $\|Ta\|_{H_2} \le C$ for all atoms $a\in H_1$. 
By linearity and definition of $S$, we have $ \|Tf\|_{H_2} \le 10 C \|f\|_{H_1}$ for all $f\in S$ and we conclude by density. \qed

\section{Boundedness of  $T$}

Let us come back to the situation where $M$ is a Riemannian manifold and $X=\r_+^*\times M$.
Before proving the     regularity of the operator $T,$ that is its boundedness on $H^1_z(X)$, let us first give some lemmas about the heat kernel $p_t(x,y)$. Indeed, using the definition of $T,$ we can write 
$$Tf(t,x)=\int_0^t \Delta e^{(t-s) \Delta} f(s,x)ds = \iint_X  p'_{t-s}(x,y) \chi _{\{t> s\}} f(s,y)ds d\mu(y),$$
where $p'_{t}(x,y)$ stands for $ \frac{\partial}{\partial t} p_{t}(x,y).$   The last equality is a formal one as the integral may not converge.

\begin{lem} (\cite{s}, Proposition 3.3) \label{lem:gub} If $M$ satisfies the doubling property $(D)$ and
the $Poincar\acute{e}$ inequality $(P)$, then the heat kernel satisfies
the Gaussian upper estimate:
$$
p_t(x,y)\leq \frac{C_1}{V(x,\sqrt{t})} \cdot e^{-c_1 \frac{{d(x,y)}^2}{t}}\eqno(G)
$$
for all $t>0$ and $x,y\in M$.

\end{lem}

There is also a lower estimate but we do not need it. 

\begin{lem}\label{lem:timept} 
(\cite{CD}, Lemma 2.5) If the heat kernel satisfies $(G)$, then the time
derivatives of $p_t$ satisfy:
\begin{equation}\label{gt}
|\frac{\partial^{k}}{\partial t^{k}}p_{t}(x,y)|
\lesssim \frac{1}{t^k V(x,\sqrt{t})}
e^{-\frac{C{d(x,y)}^2}{t}},
\end{equation}

for all $ t>0$, $x,y \in M$ and all integer $k\ge 1$.
\end{lem}





\begin{lem} \label{lem:holder} If (D) and (P) hold, for all $ t>0$, $x,y,x_0\in M$, if $d(y,x_0)\leq\sqrt{t}$, then
$$
\Big|\frac{\partial}{\partial t} p_t(x,y)- \frac{\partial}{\partial
t}p_t(x,x_0)\Big|\lesssim  \frac{1}{t V(x,\sqrt{t})}  \left(\frac{d(y,x_0)}{\sqrt{t}}\right)^{\gamma  }
e^{-\frac{C{d(x,y)}^2}{t}},
$$
where $\gamma  \in (0,1)$, $C>0$ are independent of $x,y,x_0,t$.
\end{lem}

 \noindent {\bf Proof:}  Following \cite{g} (see also \cite{r}), one can write
 $$ 
\frac{\partial}{\partial t }(p_t(x,y) - p_t(x,x_0))=\int_M \frac{\partial}{\partial t } p_{t/2}(x,z)(p_{t/2}(z,y)-p_{t/2}(z,x_0)) dz 
$$
and use the the estimate   \eqref{gt} together with the H\"older estimate of $p_t(x, y)$, as $p_t(x,y)$ is a solution of heat equation, 
 from Harnack inequality (see \cite{s}, Proposition 3.2). This is where we use (indirectly) the assumption on Poincar\'e's inequality. \qed

\begin{lem}\label{lem:CD} Assume (D) and  (G). One has for all $y \in M$ and $s>0.$
$$
\int_M |\nabla _xp_s(x,y)|d\mu(x) \lesssim  \frac{1}{\sqrt{s}}.
$$
\end{lem}

\noindent {\bf Proof:} It is proved in \cite[Lemma 2.3]{CD1} 
that there is $\gamma >0$ such that
$$
\int_M |\nabla_x p_s(x,y)|^2 e^{\gamma \frac{d^2(x,y)}{s} }d\mu(x) \lesssim \frac{1}{sV(y,\sqrt{s})},
$$for all $y \in M$ and $s>0.$ It suffices to use Cauchy-Schwarz inequality. \qed

\

Now we can prove the first part of Theorem \ref{thm:h1reg}, namely that
 $T$ is bounded from
$H^1_z(X)$ to itself. By Proposition \ref{prop:h1z=h1cw}, we prove that $T$ is bounded from $H^1_z(X)$ to $H^1(X)$. To this end, it is enough by Proposition \ref{prop:bdd} to check that $T$ maps atoms to molecules defined as above, up to a multiplicative factor.

 Suppose $a \in H^1_z(X)$ is an $(1, \infty)-$atom : there is a ball $Q= \{(t,x) \in X: \tilde{d}(t,x;t_0,x_0)< r \}$ contained in $X$ such that  supp$\,a \subset Q ,$  $\|a\|_\infty \le \nu(Q)^{-1}$ and $\int_Q a \, d\nu=0$ .
Since $T$ is bounded from $L^2(X)$ to $L^2(X)$,  
$$
 \left(\int_{4 Q \cap{X}}|Ta|^2 d\nu\right)^{\frac{1}{2}} \le C \left(\int_{{X}}|a|^2 d\nu \right)^{\frac{1}2} \leq C \nu(Q)^{-\frac{1}{2}} 
\leq c\nu (4Q\cap{X})^{-\frac{1}{2}}.
$$
Thus $c^{-1} Ta$ satisfies (2) with $j=1$ in Definition \ref{defn:molecule}.

Now to prove that $Ta$ satisfies (2) with $j\ge 2$ in Definition \ref{defn:molecule}, up to some multiplicative constant, it suffices to prove that there exists $\alpha>0$ such that for every $j\ge 2$ and 
$(t,x)\in B_j(Q)$,
$$|Ta(t,x)| \lesssim 2^{-j\alpha} \nu( 2^{j+1}Q\cap{X})^{-1}.$$
Observe that $Ta(t,x)=0$ if $t\le t_0-r^2$. Hence we assume $t>t_0-r^2$ in the sequel and we have
$$
Ta(t,x) = \int_{t_0-r^2}^t\int_{d(y,x_0)<r}p'_{t-s}(x,y) a(s,y)ds d\mu(y).
$$

 If  $(t,x)\in B_j(Q)$, then $
2^jr \leq \tilde{d}(t,x;t_0,x_0) \leq 2^{j+1}r,$ and we have the
following two cases:\\

Case $\textrm{I}$:  $ 2^jr \leq d(x,x_0) \leq
2^{j+1}r,|t-t_0| \leq (2^{j+1} r)^2 .$ 

Case $\textrm{II}$: 
$
d(x,x_0)\leq 2^{j+1}r,(2^j r)^2 \leq |t-t_0| \leq (2^{j+1}r)^2
.$\\

Case  I: we dinstiguish two subcases.

$1^{\circ}$: If $2^j r\leq d(x,x_0)\leq
2^{j+1} r$ and $|t-t_0| \leq (2r)^2,$ then we have that
$2^{j-2}r \leq d(x,y)\leq 2^{j+2}r,$ whenever $(s,y) \in Q.$ As a consequence,
\begin{equation}\label{eq1}
\begin{array}{l}

\displaystyle |Ta(t,x)|

\leq  \displaystyle \int_{t_0-r^2}^{t} \int_{d(y,x_0)<r}

\frac{1}{(t-s)V(x,\sqrt{t-s})} e^{-\frac{Cd(x,y)^2}{t-s}} |a(s,y)| ds d\mu(y)\\[0.1in]

\lesssim \displaystyle \frac{V(x_0,r)}{\nu (Q)} \int_{0}^{t-t_0+r^2} \frac{1}{s V(x,\sqrt{s})} e^{- \frac{C(2^j r)^2}{s}}ds \\[0.1in]

\lesssim  \displaystyle  \frac{V(x_0,r)}{\nu(Q)} \int_{\frac{4^j}{5}}^{+\infty} \frac{1}{uV(x,\frac{2^{j}r}{\sqrt{u}})} e^{-Cu}du \\[0.1in]

\lesssim  \displaystyle \frac{2^{j (\delta +2) }}{\nu (2^{j+1}Q)}  \int_{\frac{4^j}{5}}^{+\infty} \frac{V(x_0,r)}{uV(x,\frac{2^{j}r}{\sqrt{u}})} e^{-Cu}du \\[0.1in]

\lesssim  \displaystyle \frac{2^{j (\delta +2) }}{\nu(2^{j+1}Q)} \int_{\frac{4^j}{5}}^{+\infty} \frac{1}{u} \bigg(\frac{\sqrt{u}}{2^j}\bigg)^\delta e^{-Cu}du \\[0.1in]

\lesssim \displaystyle \frac{1}{\nu(2^{j+1}Q)}\ e^{-C'4^j} \lesssim \displaystyle \frac{1}{\nu( 2^{j+1}Q\cap X)}e^{-C''4^j}.\\

\end{array}
\end{equation}
The first inequality is obtained by Lemma \ref{lem:timept}, the second  follows from the size condition of $a$, the third  is by changing variable: $u = \frac{(2^jr)^2}{s},$ the fourth is by using doubling property (\ref{eq:doublingnu}) for  $\nu$  and the fifth  is by the doubling property of $\mu$: 
$$ \displaystyle V(x_0,r) \lesssim V(x,r) \left(1 + \frac{d(x_0,x)}{r}\right)^\delta 
\lesssim V\left(x,\frac{2^j r}{\sqrt{5u}}\right) \left(1 + \frac{2^jr}{r}\right)^\delta 
\left(\frac{\sqrt{u}}{2^j}\right)^\delta.$$

$2^{\circ}$: If $2^j r\leq d(x,x_0)\leq 2^{j+1}r$ and $ (2r)^2 \leq |t-t_0| \leq (2^{j+1}r)^2$, then
$d(x,y) \sim 2^jr$ and $t-s \sim t-t_0$ when $(s,y)\in Q$. In this case, using $\iint a(s,y)dsd\mu(y)=0$,  write $Ta(t,x)=I_1+I_2$ where
$$
I_1=  \displaystyle \int_{t_0-r^2}^{t_0+r^2} \int_{d(y,x_0)<r} (p'_{t-s}(x,y)-
p'_{t-s}(x,x_0)) a(s,y)dsd\mu(y),$$
and
$$
I_2=  \displaystyle \int_{t_0-r^2}^{t_0+r^2} \int_{d(y,x_0)<r} (p'_{t-s}(x,x_0) - p'_{t-t_0 - r^2}(x,x_0)) a(s,y) dsd\mu(y) .
$$

For $I_1$, we have the following estimate:
\begin{align*}
|I_1| &\leq \iint_Q |p'_{t-s}(x,y) - 
p'_{t-s}(x,x_0)| |a(s,y)|ds d\mu(y)\\
&
\lesssim
\displaystyle \frac{1}{\nu(Q)} \iint_Q \frac{1}{(t-s)V(x,\sqrt{t-s})}  \left(\frac{d(y,x_0)}{\sqrt{t-s}} \right)^{\gamma }
e^{- \frac{C{d(x,y)}^2}{t-s}}ds d\mu(y)\\
&
\lesssim
\displaystyle \left( \frac{r}{\sqrt{t-t_{0}}} \right)^{\gamma } \frac{1}{|t-t_{0}|V(x,\sqrt{|t-t_{0}|})}
e^{-\frac{C(2^jr)^{2}}{t-t_{0}}}\\
&
\lesssim \displaystyle \frac{1}{2^{j \gamma } \nu(2^{j+1}Q)} \left (\frac{2^jr}{\sqrt{t-t_0}} \right)^\gamma 
 \frac{(2^j r)^2 V(x_0,2^jr)}{|t-t_{0}| V(x,\sqrt{|t-t_{0}|})}
e^{-\frac{C(2^jr)^2}{t-t_0}}
\\
&
\lesssim
 \displaystyle  \frac{1}{2^{j\gamma } \nu(2^{j+1}Q)}
\\
&
\lesssim
\displaystyle  \frac{1}{2^{j \gamma } \nu( 2^{j+1}Q \cap{X})}.
\end{align*}
The second inequality is obtained by Lemma \ref{lem:holder}, and the fifth inequality is by doubling property $V(x_0,2^j r) \lesssim 
V(x,2^{j+1}r) \lesssim (\frac{2^jr}{\sqrt{t-t_0}})^{\delta } V(x,\sqrt{t-t_0})$, and that the function $ u \rightarrow  u^\beta e^{-Cu}$ is bounded on $\mathbb{R}_+$ .

Similarly, for $I_2$, we have
\begin{align*}
|I_2| &\leq \bigg| \int_{s = t_0-r^2}^{t_0+r^2} \int_{d(y,x_0)<r} \int_{w= s}^{t_0+r^2}p''_{t-w}(x,x_0)dw\ a(s,y)ds d\mu(y)\bigg|\\
&
\lesssim \displaystyle \frac{V(x_0,r)}{\nu(Q)} \int_{ s= t_0-r^2}^{t_0+r^2} \int_{ w= s}^{t_0+r^2}|p''_{t-w}(x,x_0)|dwds\\
&
\lesssim \displaystyle \frac{V(x_0,r)}{\nu(Q)} \int_{w = t_0-r^2}^{t_0+r^2} \int_{s= t_0-r^2}^{w} \frac{1}{(t-w)^{2}V(x,\sqrt{t-w})}
e^{- \frac{C{d(x,x_0)^2}}{(t-w)}}dwds \\
&
\lesssim \displaystyle \frac{V(x_0,r)}{\nu(Q)}\int_{t_0-r^2}^{t_0+r^2}(w-t_0+r^2)
\frac{1}{(t-w)^{2}V(x,\sqrt{t-w})} e^{-C \frac{(2^jr)^2}{t-w}} dw\\
&
\lesssim  \displaystyle \frac{V(x_0,r)}{\nu(Q)} \int_0^{2r^2} wdw 
\frac{1}{(t-t_{0})^{2}V(x,\sqrt{t-t_{0}})}
e^{-C \frac{4^jr^2}{t-t_0}}\\
&
\lesssim \displaystyle \frac{r^2V(x_0,r)}{\nu(Q)}\cdot \frac{1}{(2^jr)^2 V(x_0,2^jr)} \cdot
\frac{r^2(2^{j}r)^2 V(x,2^{j}r)}{(t-t_{0})^{2}V(x,\sqrt{t-t_{0}})}
e^{-C \frac{4^jr^2}{t-t_0}} \\
&
\lesssim \displaystyle  \frac{1}{ 4^j \nu(2^{j+1}Q)} \lesssim  \frac{1}{4^j \nu(2^{j+1}Q\cap X)}, 
\end{align*}
where the third inequality is obtained by Lemma \ref{lem:timept}, the fifth inequality by subsitution and the fact $t-w \sim t-t_0,$ the sixth inequality  by $V(x,2^jr) \sim V(x_0,2^jr)$, the seventh  by $ r^2V(x_0,r) \sim \nu(Q),$ and that the function $ u \rightarrow  u^\beta e^{-Cu}$ is bounded on $\mathbb{R}_+$.
\\

Case \textrm{II}: $ (2^jr)^2 \leq |t-t_0| \leq (2^{j+1}r)^2$ and $ d(x,x_0) \leq 2^{j+1}r.$
Then we have $ t-s\sim t-t_0.$ Write again $ Ta=I_1+I_2,$ where $I_1$ and $I_2$ are as in Case I. Thus
\begin{align*}
|I_1|& \leq  \displaystyle \frac 1{\nu(Q)}\int_{t_0-r^2}^{t_0+ r^2}\int_{d(y,x_0)<r}
\left (\frac{r}{\sqrt{t-s}} \right)^{\gamma } \frac{1}{(t-s)V(x,\sqrt{t-s})} e^{-\frac{C d(x,y)^2}{t-s}}ds d\mu(y)\\
&
\lesssim  \displaystyle \frac{1}{\nu(Q)} \int_{t_0-r^2}^{t_0+r^2} \int_{d(y,x_0)<r}\left(\frac{r}{\sqrt{t-s}}\right)^{\gamma }
\frac{1}{(t-s)V(x,\sqrt{t-s})}ds d\mu(y)\\
&
\lesssim
 \displaystyle \frac{1}{2^{j\gamma}} \frac{1}{(2^{j}r)^{2}V(x,2^{j}r)}
\lesssim  \displaystyle \frac{1}{2^{j\gamma } \nu(2^{j+1}Q)} \lesssim
\frac{1}{2^{j\gamma }\nu(2^{j+1}Q\cap X)},
\end{align*}
where the first inequality follows from Lemma \ref{lem:holder}, the third inequality follows from the fact that $t-s \sim (2^jr)^2,$ the fourth from $V(x,2^jr) \sim V(x_0,2^jr)$.

For $I_2,$ we have 
\begin{align*}
|I_2| &\leq \displaystyle \int_{s= t_0-r^2}^{t_0+r^2} \int_{d(y,x_0)<r}\bigg\{ \int_{\theta = s}^{t_0+r^2}|{p''_{t-\theta}}(x,x_0)| d\theta\bigg\}
|a(s,y)|ds d\mu(y) \\
&
\lesssim  \displaystyle \frac{1}{\nu(Q)} \int_{t_0-r^2}^{t_0+r^2} \int_{d(y,x_0)<r} \int_{s}^{t_0+r^2}
\frac{1}{(t-\theta)^{2}V(x,\sqrt{t-\theta})} d\theta
ds d\mu(y)\\
&
\lesssim  \displaystyle \frac{r^2}{(2^{j}r)^{4}V(x_{0},2^{j}r)}
\lesssim \frac{1}{4^j \nu( 2^{j+1}Q)}
\lesssim \frac{1}{4^j \nu(2^{j+1}Q \cap X)},
\end{align*}
where the second inequality is obtained by Lemma \ref{lem:timept}, the third one  by $ t-\theta \sim (2^jr)^2$ and $V(x,2^jr) \sim V(x_0,2^jr)$.

Thus we proved that $Ta$ satisfies (2) with  $\alpha=\gamma$ in Definition \ref{defn:molecule}.

We claim that $\int_X Ta\,  d\nu = 0$. Assuming this claim for a moment, we conclude from Proposition \ref{lem:molecule}  that $Ta \in H^1(X)$ with $\|Ta\|_{H^1(X)} \le C$ uniformly with respect to $a$. \hfill${\Box}$

\begin{lem}\label{lem:meanvalue} $\int_X Ta \, d\nu =0$ when $a \in H_z^1(X)$ is an $(1,\infty)$-atom.
\end{lem}

\noindent {\bf Proof:} We give two proofs, each having advantages. 
The first one is as follows. Let $u$ be the solution (since  $a \in L^2(X)$) of \eqref{eq0} with data  $a$. 
Then $Ta=\Delta u =  u' - a$. Furthermore $a \in L^1(X)$ and $Ta, u' \in L^1(X)$ as we just showed. Fix $h$ a bounded function with compact support on $M$.  
Then integrate the equation against $h$ on $X$. Note that  $u\in C([0, \infty), L^2(M))$     with $u(0,\cdot)=0$ and $\lim_{t\to \infty} u(t, \cdot) =0$ in $L^2(M)$
(since $a$ is 0 for $t$ larger than some $\tau$, hence $u(t, \cdot) = e^{(t-\tau)\Delta} u(\tau, \cdot)$ fot $t> \tau$) and $u' \in L^2(\r_+^*, L^2(M))$. Integrating by parts, one finds
\begin{equation*}
\iint_X u'(t,x) h(x) \, dtd\mu(x)=    \int_0^\infty < u'(t,\cdot), h>\, dt = 0,
\end{equation*}hence,  
 \begin{equation*}
\iint_X Ta(t,x) h(x) \, dtd\mu(x)=  -  \iint_X a(t,x) h(x) \, dtd\mu(x).
\end{equation*}
Since $a$ and $Ta$ both belong to $L^1(X)$, letting $h$ converge to $1$ everywhere  with $\|h\|_\infty =1$ implies
  \begin{equation*}
\iint_X Ta(t,x)  \, dtd\mu(x)=  -  \iint_X a(t,x)  \, dtd\mu(x)=0.
\end{equation*}

The second argument is more precise: by the previous calculations, we know that
$$
\iint_X |Ta(t,x)| dtd\mu(x) \lesssim
\sum_{j=1}^{\infty} 2^{-j \alpha} < +\infty,$$
 thus
$Ta \in L^1(X)$ and, therefore, for a.e.
$t>0,~~ Ta(t,\cdot) \in L^1(M)$. It is enough to show that for such $t$'s, $\int_M Ta(t,x)\, d\mu(x)=0$.

 Now let $\varphi \in
C_0^{\infty}([0,+\infty)),$ with $\varphi=1$ on $[0,1]$, $\varphi=0$ on $[2,+\infty)$ and $0\leq \varphi \leq1$. Fix $x_0\in M$ and set $\varphi_R(x)=\varphi(R^{-1} \eta(x))$ for $R>0$ where $\eta\in C^\infty(M\setminus \{x_0\})$ with $\eta(x)\sim d(x_0,x)$ and $\nabla \eta$ bounded. 
Recall that $Ta=\Delta u$ where $u(t)=\int_0^t e^{(t-s)\Delta}a(s)\, ds$. Hence
\begin{equation}\label{eq98}
\int_{M} Ta(t,x) \varphi_R(x) d\mu(x)= -\int_{M} \nabla
u(t,x) \cdot \nabla \varphi_R(x) d\mu(x).
\end{equation}
Since
$ \|\nabla \varphi_R \|_{\infty} \leq \frac{C}{R},$
we obtain
\begin{equation}\label{eq99}
\left |\int_{M} Ta(t,x) \varphi_R(x) d\mu(x) \right| 
\lesssim \frac{1}{R} \int_M |\nabla u(t,x)|d\mu(x).
\end{equation}
It follows immediately from Lemma \ref{lem:CD} that 
\begin{equation*}
\|\nabla e^{ s \Delta }g\|_{L^1(M)} \lesssim \frac{1}{\sqrt{s}} \|g\|_{L^1(M)}, 
\end{equation*}
for all $g \in L^1(M)$ and $s>0$, therefore 

\begin{equation} \label{eq100}
\int_M |\nabla u(t,x)|d\mu(x) \leq  \int_0^t \frac{C}{\sqrt{t-s}} \|a(s,\cdot)\|_{L^1(M)} ds 
\leq 2C \sqrt t\,  \|a\|_{\infty}r^{-2},
\end{equation}
where the ball supporting $a$  has radius $r$.
Thus by (\ref{eq98}), (\ref{eq99}), (\ref{eq100}), we have
$$
\lim_{R \rightarrow \infty} \left|\int_M Ta(t,x) \varphi_R(x) d\mu(x) \right|=0.$$
Hence, as $t>0$ is such that $Ta(t,\cdot) \in L^1(M)$, dominated convergence yields 
$$
\int_{M} Ta(t,x) d\mu(x) =\lim_{R \rightarrow \infty} \int_{M} Ta(t,x) \varphi_R(x) d\mu(x) =  0.$$

 \hfill{$\Box$}     
 
 \begin{rem} The second argument shows that for almost all $t>0$, $Ta(t, \cdot)$ is in $L^1(M)$ and has a vanishing integral. Actually, an examination of the proof (left to the reader)  shows that for almost all $t>0$, $Ta(t,\cdot)$ is  up to a constant depending on $t, t_0, r$ a molecule 
 in $M$, hence belongs to $H^1_{CW}(M)$. Calling $c_{t_0,r}(t)$ its norm in this space,  one can check that 
 $\int_0^\infty c_{t_0,r}(t) \, dt \lesssim 1$ uniformly with respect to $t_0,r$, that is uniformly with respect to the atom $a$. Thus, the maximal regularity map $T$ extends boundedly from $H^1_z(X)$ into $L^1(\r_+^*, H^1(M))$.  It is easy to show that $H^1_z(\r_+^*\times M)$ and $L^1(\r_+^*, H^1(M))$ are not comparable spaces. This extra regularity of the map $T$ is thus not shared by the
 map $f\mapsto u'$.  
 
 \end{rem}

We conclude this section by observing that the source space $H^1_z(X)$ is best possible.

\begin{prop} If $M=\r^n$, 
$T$ is not bounded from $H^1_r(X)$ to $L^1(X).$  
\end{prop}

\noindent {\bf Proof:} Let $f(t,x) = 1$ if $0<t<1$ and $|x|\leq  1$, and $0$ otherwise. We have seen  that 
$f \in H^1_r$. We claim that $Tf \notin L^1(X)$. Note that $X=\r_+^*\times \r^n$.

Recall that if $p_t(z)$ denotes the heat kernel, then  the time derivative is given by
$$
p'_t(z) = \frac 1 {(4\pi t)^{n/2} } \left( \frac{|z|^2 }{4t^2} - \frac n {2t}\right) e^{- \frac {|z|^2}{4t}},
$$
hence 
$$
 \frac{|z|^2 }{4t}  \le \frac  n  4 \Longrightarrow p'_t(z) \le -  \frac {n} {4t (4\pi t)^{n/2}}.  $$
Thus, if $|x| \le \frac 1 2 {(nt)}^{1/2}$, $|y|\le 1$, $t\ge \sup (4/n, 2)=\tau$ and $0<s<1$ we have
$$
 p'_{t-s}(x-y) \le -  \frac {n} {4t (4\pi t)^{n/2}}  .
 $$
 It follows that there is a constant $c_n>0$ such that
\begin{align*}
\iint |Tf(t,x)| dtdx) &\ge \int_\tau^\infty  \int_{|x| \le \frac 1 2 {(nt)}^{1/2}} \left|\int_{|y|\le 1} \int_0^1 p'_{t-s}(x-y)dsdy\right|dtdx 
\\
& \ge c_n  \int_\tau^\infty \frac {dt}t
\\
&=+ \infty.
\end{align*}

\section{Boundedness of $T^{*}$}

In this section, we  study the boundedness of $T^{*}$ given by (\ref{eq:T*}).  If we let 
$$v(t,x) = \int_t^{+\infty}
\big[ e^{(s-t)\Delta} f(s,\cdot)\big](x) ds, t>0,$$
 then $\Delta v=T^*f$ and
 satisfies: 
$$\frac{\partial v}{\partial t}(t,x) + \Delta v(t,x) = -f(t,x).$$

We now  prove the part of Theorem \ref{thm:h1reg} concerning the boundedness of $T^*$. 

 \noindent {\bf Proof:} As before, it is enough to prove the boundedness from $H^1_r(X)$ into $H^1(X)$ and for that purpose that   $H^1_r$-atoms are mapped to $H^1$-molecules of Definition \ref{defn:molecule} up to a multiplicative constant. There are two kinds of atoms in
$H_r^1(X)$: type $(a)$ atoms and type $(b)$ atoms.

Suppose $a$ is a type $(a)$ atom, with supp$\, a \subseteq
Q=\{(t,x):\widetilde{d}(t,x;t_0,x_0)<r \}$ and $4Q\subseteq X$. To see that $T^*a$ is a molecule we  prove (1) and (2) in Definition \ref{defn:molecule}. We begin with (2). First, since $T$ is
$L^2(X)$-bounded, so is $T^{*}$  and
\begin{equation}\label{eq2}
\|T^{*}a\|_{L^2(4Q \cap
X)}\lesssim{\nu(Q\cap X)^{-\frac{1}{2}}} \lesssim{\nu(4Q\cap X)^{-\frac{1}{2}}}.
\end{equation}

 Next, we show that  there exits $\alpha > 0$ such that for each $j\ge 2$ and  $(t,x)\in B_j(Q)$,
\begin{equation}\label{eq3}
|T^{*}a(t,x)| \lesssim {2^{-j\alpha}} {\nu( 2^{j+1}Q\cap
X)}^{-1}.
\end{equation}
We observe that $T^*a(t,x)=0$ if $t>t_0+r^2$. Hence we assume $t\le t_0+r^2$.
For $(t,x)\in B_j(Q)$, we have that
$2^j r\leq \widetilde{d}(t,x;t_0,x_0) \leq 2^{j+1}r$.
So we distinguish in addition two cases:\\

 Case $\textrm{I}$:  $ 2^j r\leq
d(x,x_0)\leq 2^{j+1}r, |t-t_0| \leq (2^{j+1}r)^2\ .$

Case $\textrm{II}$:   
 $
(2^jr)^2 \leq |t-t_0| \leq (2^{j+1}r)^2, d(x,x_0) \leq 2^{j}r\
.$\\

Case $\textrm{I}$: 
 $ 2^j r\leq d(x,x_0) \leq 2^{j+1}r$ and $|t-t_0|\leq (2^{j+1}r)^2.$ If, in addition $|t-t_0|\le (2r)^2$ then
\begin{align*}
|T^{*}a(t,x)| &=  \left| \int_t^{+\infty} \int_{d(y,x_0)<r} 
p'_{s-t}(x,y)a(s,y)ds d\mu(y)\right|
\\
&
\lesssim  \frac{1}{\nu(Q)} \int_t^{t_0+r^2} \int_{d(y,x_0) < r}
\frac{1}{(s-t)V(x,\sqrt{s-t})} e^{- \frac{C {d(y,x)}^2}{s-t}}ds d\mu(y)\\
&
\lesssim  \frac{V(x_0,r)}{\nu(Q)} \int_t^{t_0 + r^2} \frac{1}{(s-t)V(x,\sqrt{s-t})}
e^{-\frac{C4^jr^2}{s-t}}ds \lesssim
e^{-C4^j} \frac{1}{\nu(2^jQ\cap X)},
\end{align*}
where  we used Lemma \ref{lem:timept} and that $a$ is supported in $Q.$ The last inequality is obtained repeating  the proof of (\ref{eq1}).

If $ t \leq t_0- (2r)^2$, as $t \le t_0+r^2$, we have $t\le t_0-4r^2$, then
\begin{align*}
|T^{*}a(t,x)| &\leq  \frac{1}{r^2} \int_{t_0-r^2}^{t_0+r^2} \frac{1}{(s-t)V(x,\sqrt{s-t})}
e^{- \frac{C4^j r^2}{s-t}}ds \\&
\leq
  \frac{V(x_0,r)}{\nu(Q)} \int_{\frac{4^jr^2}{t_0+ r^2 -t}}^{\frac{4^jr^2}{t_0-r^2-t}} \frac{e^{-Cu}du}{uV(x,\frac{2^{j}r}{\sqrt{u}})}
\\
&
\lesssim   \frac{V(x_0,r)}{\nu(Q)} \int_{\frac{4^j}{2}}^{+\infty} \frac{e^{-Cu}du}{uV(x,\frac{2^jr}{\sqrt{u}})}
\\
&
\lesssim e^{- c4^j} \frac{1}{\nu(2^{j+1}Q \cap X)},
\end{align*}
where the second inequality is obtained by substitution, the third  from the fact that $t_0-t + r^2 \geq r^2,$
 and the last one  repeating the proof of (\ref{eq1}). \\

Case $\textrm{II}$: 
$(2^jr)^2 \leq |t-t_0| \leq (2^{j+1}r)^2, d(x,x_0) \leq 2^{j}r .
$ In this case, using the vanishing moment condition of $a$, we can write $T^{*}a(t,x)$ as follows:
\[
\begin{array}{l}
T^{*}a(t,x) = \displaystyle \int_{t_0-r^2}^{t_0+r^2} \int_{d(y,x_0)\leq r} 
p'_{s-t}(x,y)a(s,y)ds d\mu(y)\\[0.1in]

=  \displaystyle \int_{t_0-r^2}^{t_0+r^2} \int_{d(y,x_0)\leq r}(p'_{s-t}(x,y)- p'_{s-t}(x,x_0)) a(s,y)ds d\mu(y)  \\[0.1in]

 \displaystyle  +  \int_{t_0-r^2}^{t_0+ r^2} \int_{d(y,x_0)\leq r} (p'_{s-t}(x,x_0) - p'_{r^2+t_0-t}(x,x_0)) a(s,y)ds d\mu(y)
=J_1+ J_2.
\end{array}
\]
Following the  proof of (\ref{eq1}), we have the estimation of $J_1$ :
\[
\begin{array}{l}
|J_1| \lesssim \displaystyle \frac{1}{\nu (Q)} \int_{t_0-r^2}^{t_0+ r^2} \int_{d(y,x_0)\leq r}\frac{1}{(s-t)V(x,\sqrt{s-t})}
e^{- \frac{C{d(y,x)}^2}{s-t}}ds d\mu(y)
\lesssim  \frac{1}{ 2^{j\gamma}\nu(2^{j+1}Q \cap X)},
\end{array}
\]
For $J_2,$ we have the estimation:
\[
\begin{array}{l}

|J_2| \leq \displaystyle \int_{t_0-r^2}^{t_0+r^2} \int_{d(y,x_0)\leq r} \int_{\theta = s}^{t_0+r^2}|p''_{\theta - t}(x,x_0)|d\theta
|a(s,y)|ds d\mu(y)\\[0.1in]

\lesssim \displaystyle \frac{1}{\nu(Q)} \int_{s = t_0-r^2}^{t_0+ r^2} \int_{d(y,x_0)\leq r} \int_{\theta = s}^{t_0+r^2}

\frac{d\theta dsdy}{(\theta - t)^{2}V(x,\sqrt{\theta -t})}
\lesssim \frac{1}{4^j}\cdot\frac{1}{\nu(2^jQ\cap X)},
\end{array}
\]
where we use Lemma \ref{lem:timept} and the fact that $ \theta  - t \sim (2^jr)^2$ to get the last inequality.
 
Now the argument to show $\int_X T^{*}a(t,x)dtd\mu(x)=0$ is the same as in Lemma \ref{lem:meanvalue}. This proves that $T^{*}a$ is a molecule.

Now we prove that for type $(b)$-atoms $b$, $T^{*}b$ is a molecule.
Suppose that $b$ is a type $(b)$ atom, supported in
$Q=\{(t,x) \in X, \widetilde{d} (t,x; t_0,x_0)< r \}$ with $2Q\subseteq X$ but $4Q \nsubseteq X$, and
$\|b\|_{\infty} \leq \frac{1}{\nu(Q)}$.

As for type $(a)$ atoms, we have 
\begin{equation}\label{eq3}
\|T^{*} b\|_{L^2(4Q\cap X)} \lesssim \nu(Q)^{-\frac{1}{2}} \lesssim \nu(4Q\cap X)^{\frac 1 2}.
\end{equation}
Next, we establish that for all $j\ge 2$ and $(t,x) \in B_j(Q)$,
\begin{equation}\label{eq4}
|T^{*}b(t,x)| \lesssim \nu( 2^jQ \cap X)^{- \frac{1}{2}}e^{-C4^j}.
\end{equation}

Fix $j \geq 2$ and suppose $(t,x) \in
B_j(Q)$. Observe that $T^*b(t,x)=0$ if $t>t_0+r^2$. As $t_0\le (4r)^2$, this implies that $|t-t_0| \le 16r^2$. Thus the condition $2^j r\leq \widetilde{d}(t,x;t_0,x_0) \leq 2^{j+1}r$ implies  $d(x,x_0) \sim 2^j r$ for $j\ge 5$ (we skip the details for $j \|e 4$ which are not difficult). Like the proof of (\ref{eq1}), we have 
\begin{align*}
|T^{*} b(t,x)| &\leq  \int_t^{t_0+ r^2} \int_{d(y,x_0) \leq r}
|p'_{s-t}(x,y)||a(s,y)|ds d\mu(y)\\
&
\lesssim  \frac{1}{\nu( Q)} \int_t^{t_0+ r^2} \int_{d(y,x_0)\leq r}
\frac{1}{(s-t)V(x,\sqrt{s-t})}e^{- \frac{C {d(x,y)^2}}{s-t}}ds d\mu(y)\\
&
\lesssim   \frac{V(x_0,r)}{\nu(Q)} \int_t^{t_0+ r^2} \frac{1}{(s-t)V(x,\sqrt{s-t})}
e^{-\frac{C4^j}{s-t}}ds\\
& \lesssim e^{-C4^j} \frac{1}{\nu(2^jQ)}.
\end{align*}
The second argument in the proof of Lemma \ref{lem:meanvalue}  implies  that $\iint_X T^{*} b\, d\nu = 0$. Hence, $T^{*}b$ is, up to a multiplicative constant, a molecule  in $H^1(X)$.
The theorem is proved.  \hfill$\Box$

\begin{rem} As for $T$, one can check that the proof yields that $T^*$ is bounded
from $H^1_r(X)$ into $L^1(\r_+^*, H^1(M))$.
\end{rem}

\begin{rem} If $X=\r_+^*\times \r^n$, $T^*$ is not bounded on $L^1(X)$. Indeed, this would mean that 
$\iint_{X} \chi_{s>t} |p'_{s-t}(x-y)|\, dsdy $ is bounded uniformly with respect to $(t,x) \in X$. But this integral equals  $c \int_0^\infty \frac {dt}t=\infty$ as $ c=\int_{\r^n} |p'_1(x)|\, dx >0$. 
\end{rem}

\section{Maximal $L^q$    regularity on $L^p$}

The well-known  maximal   $L^q$    regularity of  the Cauchy problem in our setting (see \cite{DS}) states as follows 

\begin{thm} Assume that $M$ satisfies (D) and (P).
 There is   maximal  $L^q$    regularity property for the Cauchy problem (\ref {eq0}):
$$\|u'\|_{L^q(\mathbb{R}_+, Y)} + \|T f\|_{L^q(\mathbb{R}_+, Y)} \lesssim \|f\|_{L^q(\mathbb{R}_+, Y)},$$
where 
$ Y = L^p(M),$ $1<p,q<\infty$.
\end{thm}

Let us see a proof as a corollary of Theorem \ref{thm:h1reg} (as said in the introduction,  \cite{CD} gives a proof under weaker assumptions).

 \noindent {\bf Proof:} It suffices to prove that $T$ is $L^p(X)-$bounded, with $1< p < \infty.$ To see this, we construct an operator $\widetilde T$ as follows:

$$\left\{\begin{array}{ll}
\widetilde T f(t,x)= T( \chi_X f)(t, x), &  t>0,\\
\widetilde T f(t,x)=0, & t \leq  0.
\end{array}\right.
$$
Then by Theorem \ref{thm:h1reg}, we know that $T$ is bounded from $H^1_z (X)$ to $H^1_z(X),$ therefore $\widetilde T$ is bounded from $H^1 (\tm)$ to $H^1 (\tm).$ On the other hand $T$ is $L^2(X)-$bounded, so $\widetilde T$ is $L^2(\tm)-$bounded too. Then using the interpolation theorem in \cite{cw}, we obtain that $\widetilde T$ is $L^p(\tm)-$bounded, where $1<p\leq 2$. Therefore $T$ is $L^p(X)-$bounded, $1<p\leq 2.$ 

Similarly, we can repeat the argument for the operator $\widetilde T^*$ defined by
$$\left\{\begin{array}{ll}
\widetilde T^* f(t,x)= T^*( \chi_X f)(t, x), &  t>0,\\
\widetilde T^* f(t,x)=0, & t \leq  0.
\end{array}\right.
$$
and obtain that $T^*$ is $L^p(X)-$bounded, $1<p\leq 2.$ 

By duality, we obtain that $T$ is $L^p(X)-$bounded, $2<p <\infty.$ 
\hfill$\Box$

\section{Sublaplacians of connected Lie groups of polynomial growth}

First let us recall some properties of connected Lie groups of polynomial growth. For further informations, see \cite{vsc}, \cite{al}, and the references given therein.

Let $G$ be a connected Lie group. Assume $G$ is unimodular and and  fix a (left-invariant) Haar measure $\mu$ on $G$. If $A$ is a measurable subset of $G$, we will denote the measure of $A$ by $|A|.$

Let $X_1, \cdots,X_n$ be left-invariant vector fields on $G$ which satisfy the H\"{o}rmander's condition, i.e. they generate, together with their successive Lie brackets $[X_{i_1},[\cdots, \\
X_{i_\alpha }] \cdots],$ the Lie algebra of $G$. 

Associated with $X_1, \cdots,X_n,$ in a canonical way, there is a  control distance $\rho,$ which is left-invariant and compatible with the topology on $G$. 

If we put $V(r) = |B(e,r)|$ for $r>0$, where $B(e,r) = \{x \in G: \rho(x,e)< r\}$,  there exists $d\in \mathbb{N},$ which depends on the vector fields,  such that
$$ V(r) \sim r^d, ~~\mbox {for} ~~ 0<r\le 1.$$
We assume that $G$ has polynomial volume growth, namely that 
$$ V(r) \sim r^D, ~~\mbox {for} ~~1 < r <+\infty.$$
Otherwise, by a theorem of Guivarc'h, $G$ would have exponential growth.
By direct caculation, we have the "doubling property" : $V(s r) \lesssim  s ^{max \{d, D\}}  V(r),$ where $s>1.$


For example every connected nilpotent Lie group has polynomial volume growth.\\

On connected Lie groups of polynomial growth, associated with the  sublaplacian $\Delta =\sum_{j=1}^n X_j^2$ and corresponding heat kernel, Poincar\'e
inequalities also hold (see \cite{v}, \cite{s}). The heat kernel satisfies the Gaussian upper estimate $(G)$ (see \cite{vsc}, Theorem VIII.2.9). In this setting, the Harnack inequality holds too, see \cite{al}, Theorem 3.1. Gaussian estimates (\ref{gt}) of time derivatives of the heat kernel are obtained in \cite{vsc}, Theorem VIII.2.4.
Next, as Harnack inequality holds,   Lemma \ref{lem:holder} is valid.

Eventually Lemma \ref{lem:CD} is easy. The simple proof is a consequence  of the following Gaussian estimates of first order space derivatives of the heat kernel.

\begin{lem}(\cite{vsc}, Theorem VIII.2.7) \label{lem:gs}
If $G$ is a connected Lie group with polynomial growth, then for $i\in \{1, \cdots, n\},$ we have
$$
|X_i p_t(x)| \lesssim t^{-1/2} V(\sqrt{t})^{-1} e^{-c \frac{{\rho(x,e)}^2}{t}},$$
for all $t>0$ and $x \in G.$

\end{lem}




So far we have recovered  all the basic tools used in Section 4. The  ``parabolic'' Hardy spaces on $N=\r \times G$ and  the $z$ and $r$ Hardy spaces on $X=\r_+^*\times G$ can be defined analogously. Consider the Cauchy problem associated with $\Delta$ and define the operators $T$ and $T^*$ as before.   

\begin{thm} Let $G$ be a connected Lie group with polynomial growth. Then $T$ is bounded on  
$H^1_z(X)$  and $T^*$ is bounded from $H^1_r(X)$ to $H^1_z(X)$.
\end{thm}

The proof is the same as before. Details are left to the reader.\

\section{Second order operators}

Let $A\in
L^{\infty}(\reel^n,M_{n}(\complex))$. Assume
that $A$ is uniformly elliptic, which means that there exists
$\delta>0$ such that, for almost all $x\in \reel^n$ and all $\xi\in
\complex^n$,
\[
\mbox{Re }A(x)\xi. \overline{\xi} \geq \delta \left\vert
\xi\right\vert^{2}.
\]
Let $\Omega$ be an open subset of $\r^n$. If $V(\Omega)$ is a closed subset of $W^{1,2}(\Omega)$
containing $W^{1,2}_{0}(\Omega)$, there is a unique operator
$L$ which is the maximal accretive operator
associated with the accretive sesquilinear form
\[
Q(f,g)=\int_{\Omega}A(x)\nabla f(x).\overline{\nabla g(x)}dx
\]
for all $(f,g)\in V\times V$. We are interested in the Dirichlet
boundary condition ($V=W^{1,2}_{0}(\Omega)$) and the Neumann boundary
condition ($V=W^{1,2}(\Omega)$). 

\medskip

Say that $L$ satisfies $(G)$ is the following three conditions hold:

The kernel of $e^{-tL}$, denoted by $K_t(x,y)$, is a measurable
function on $\Omega\times\Omega$ and there exist
$C_G,\alpha>0$ such that, for all $0<t<+\infty$ and almost every
$x,y\in \Omega$,
\begin{equation} \label{Gauss}
\left\vert K_t(x,y)\right\vert \leq \frac {C_G}{t^{n/2}}e^{-\alpha
\frac{\left\vert x-y\right\vert^2}t}.
\end{equation}

For all $y\in \Omega$ and all $0<t<+\infty$, the function $x\mapsto
K_t(x,y)$ is H\"older continuous in $\Omega$ and there
exist $C_H,\mu>0$ such that, for all $0<t<+\infty$ and all
$x,x^{\prime},y\in \Omega$,
\begin{equation} \label{Holder1}
\left\vert K_t(x,y)-K_t(x^{\prime},y)\right\vert \leq \frac
{C_H}{t^{n/2}} \frac{\left\vert
x-x^{\prime}\right\vert^{\mu}}{t^{\mu/2}}.
\end{equation}

For all $x\in \Omega$ and all $0<t<+\infty$, the function $y\mapsto
K_t(x,y)$ is H\"older continuous in $\Omega$ and there
exist $C_H,\mu>0$ such that, for all $0<t<+\infty$ and all
$y,y^{\prime},x\in \Omega$,
\begin{equation} \label{Holder2}
\left\vert K_t(x,y)-K_t(x,y^{\prime})\right\vert \leq \frac
{C_H}{t^{n/2}} \frac{\left\vert
y-y^{\prime}\right\vert^{\mu}}{t^{\mu/2}}.
\end{equation}

We assume from now on that $\Omega$ is a strongly Lipschitz domain (see \cite{ar} for a definition). Hence $\Omega$ equipped with Euclidean distance and Lebesgue measure is a space of homogeneous type. We set $N=\r\times \Omega$, $X=\r_+^*\times \Omega$ and we can define the parabolic  Hardy spaces   $H^1(N)$, $H^1_r(X)$ and $H^1_z(X)=H^1(X)$ following 
Section  \ref{sec:hardy}. Note that the subscripts $r$ and $z$ are with respect to extension with respect to the variable $t$ across $t=0$,  the spatial  variable remaining in $\Omega$. 

Consider the Cauchy problem associated with $-L$ on $X$ and call again $T$ the map of maximal regularity and $T^*$ its adjoint in $L^2(X)$. The result is as follows.

\begin{thm} Suppose $\Omega$ is a strongly Lipschitz domain and that $L$ satisfies condition 
(G).  
\begin{enumerate} 
\item In the case of Dirichlet boundary condition, $T$ is bounded from $H^1_z(X)$ into $H^1_r(X)$ and $T^*$ from $H^1_r(X)$ into $H^1_r(X)$. 
\item In the case of Neumann boundary condition, $T$ is bounded from $H^1_z(X)$ into $H^1_z(X)$ and $T^*$ from $H^1_r(X)$ into $H^1_z(X)$. 
\end{enumerate}
\end{thm}

For example, this theorem applies when $A$ has real-valued entries. 

Let us sketch the proof as it is essentially similar to the ones in the case of manifolds.

First, the four lemmas on $p_t(x,y)$ have analogs for $K_t(x,y)$. The upper bound in (G) correspond to the  Gaussian upper bound of  Lemma \ref{lem:gub}. By analyticity of the semigroup generated by $-L$, it  easily implies  Lemma \ref{lem:timept} and Lemma \ref{lem:holder}. Finally, the estimate on the gradient of Lemma \ref{lem:CD}  is a consequence of  \cite[Proposition A.4]{ar}. 

  This yields that the action of $T$ and $T^*$ on appropriate atoms are  functions with same estimates as molecules. 

For both boundary conditions the source spaces are the same.  What changes are the target spaces, let us explain why. In the case of  Neumann boundary conditions, one has $\int_\Omega \partial_t K_t(x,y)\, dx=0$ so one expects $Ta$ and $T^*a$ to have mean value zero. Indeed, this is obtained by mimicking  the second argument for Lemma  \ref{lem:meanvalue}. In the case of Dirichlet boundary condition, the first argument  for Lemma  \ref{lem:meanvalue} applies. Hence, when the corresponding atom has mean value zero, then the same is true for its image under $T$ or $T^*$. One concludes with Proposition  \ref{lem:molecule}.
Further details are left to the reader.

\begin{rem} One could also consider elliptic operators on manifolds with (D) and (P)   or subelliptic oprerators on connected Lie groups of polynomial growth with coefficients and appropriate estimates, and develop a similar theory. This is left to the reader. 
\end{rem} 
\begin{rem} All results have local analogs under local doubling, local Poincar\'e and local Gaussian bounds in which $t$ is restricted to a finite interval $(0, \tau)$. In this case, $\r_+^*$ and $\r$ are  replaced  in $X$ and $N$ by $(0, \tau)$ and $(-\tau, \tau)$ and the Hardy spaces by their local analogs in the sense of Goldberg (where one relaxes the mean value zero property for atoms or molecules associated to large balls).    This is left to the interested reader. 
\end{rem}

\

{\bf  Acknowledgements:}  The authors thank Thierry Coulhon for 
indicating useful bibliographical references and Emmanuel Russ for reading carefully the manuscript and suggesting some improvements. The third author would like to express her thanks to professor Lizhong Peng and professor Heping Liu for their encouragements, and thanks the Department of Mathematics of theUniversity Paris-Sud  for its hospitality.

\end{document}